\documentclass[10pt,a4paper]{article}

\usepackage{booktabs}
\usepackage[utf8]{inputenc}
\usepackage[english]{babel}
\usepackage{amsmath}
\usepackage{amsfonts}
\usepackage{amssymb}
\usepackage{makeidx}
\usepackage{graphicx}
\usepackage{lmodern}
\usepackage{kpfonts}
\usepackage{dsfont}
\usepackage{cancel}
\usepackage{amsthm} 
\usepackage{hyperref}
\usepackage[left=1.5cm,right=1.5cm,top=2cm,bottom=2cm]{geometry}

\usepackage{subcaption}
\usepackage{multirow}
\usepackage{float}
\usepackage{url}
\usepackage{multicol}
\usepackage{mathtools, cuted}
\usepackage{lipsum}
\usepackage{fancyhdr}
\newtheorem{theorem}{ Theorem}[section]
\newtheorem{definition}[theorem]{Definition}
\newtheorem{proposition}[theorem]{Proposition}
\newtheorem{example}[theorem]{Example}
\providecommand{\nset}[1]{
\mathbb{#1}
}
\providecommand{\set}[1]{
\left\{#1\right\}
}
\providecommand{\com}[1]{``#1"}

\providecommand{\ifr}[5]{
{}^{#1}_{#2}{#3}_{#4}^{#5}
}
\providecommand{\gam}[1]{
\Gamma\left(#1 \right)
}
\providecommand{\re}[1]{
Re\left(#1 \right)
}
\providecommand{\norm}[1]{
\left\lVert #1 \right\rVert
}

\providecommand{\abs}[1]{
\left\lvert #1 \right\rvert
}

\providecommand{\ds}[1]{
\displaystyle #1
}

\providecommand{\der}[3]{
\dfrac{{#1}^{#3}}{{#1}{#2}^{#3}}
}
\providecommand{\dt}[1]{
\mbox{det}\left(#1\right)
}
\newcommand{\mref}[2]{\textbf{#1 #2}}

\usepackage{enumitem} 
\setlist[itemize]{noitemsep} 

\usepackage{abstract} 

\usepackage{titlesec} 
\titleformat{\section}[block]{\large\scshape\centering}{\thesection.}{1em}{} 
\titleformat{\subsection}[block]{\large\scshape\centering}{\thesubsection.}{1em}{} 
\usepackage{titling} 

\pretitle{\begin{center}\Huge\bfseries} 
\posttitle{\end{center}} 
\title{Fractional Newton-Raphson Method and Some Variants for the Solution of Non-linear Systems} 
\author{
\textsc{A. Torres-Hernandez } \\ 
\normalsize  Department of Physics - UNAM \\ 
\normalsize  \href{mailto:anthony.torres@ciencias.unam.mx}{anthony.torres@ciencias.unam.mx}
\and  
\textsc{F. Brambila-Paz }\\ 
\normalsize Department of Mathematics - UNAM \\
\normalsize \href{mailto:fernandobrambila@gmail.com}{fernandobrambila@gmail.com} 
}
\date{} 
\renewcommand{\maketitlehookd}{%

\begin{abstract}

We present some novel numerical methods valid for one and several variables, which using the fractional derivative, allow to find solutions for some non-linear systems in the complex space using real initial conditions. The origin of these methods is the fractional Newton-Raphson method (see \cite{fernando2017fractional,brambila2018fractional}), but unlike the latter, the orders proposed here for the fractional derivatives are functions. In the first method, a function is used to guarantee an order of convergence (at least) quadratic, and in the others, a function is used to avoid the discontinuity that is generated when the fractional derivative of the constants is used, and with this, it is possible that the methods have at most an order of convergence (at least) linear.

\textbf{Keywords:} Iteration Function, Order of Convergence, Fractional Derivative.
\end{abstract}
}

\begin{document}

\maketitle


\section{Introduction}

When starting to study the fractional calculus, the first difficulty is that, when wanting to solve a problem related to physical units, such as determining the velocity of a particle, the fractional derivative seems to make no sense, this is due to the presence of physical units such as meter and second raised to non-integer exponents, opposite to what happens with operators of integer order. The second difficulty, which is a recurring topic of debate in the study of fractional calculus, is to know what is the order \com{optimal} $\alpha$ that should be used when the goal is to solve a problem related to fractional operators. To face these difficulties, in the first case, it is common to dimensionless any equation in which non-integer operators are involved, while for the second case different orders $\alpha$ are used in fractional operators to solve some problem, and then it is chosen the order $\alpha$ that provides the \com{best solution} based about an established criteria.

Based on the two previous difficulties, arises the idea of looking for applications with dimensionless nature such that the need to use multiple orders $\alpha$ can be exploited in some way. The aforementioned led to the study of Newton-Raphson method and a particular problem related to the search for roots in the complex space for polynomials: if it is needed to find a complex root of a polynomial using Newton-Raphson method, it is necessary to provide a complex initial condition $x_0$ and, if the right conditions are selected, this leads to a complex solution, but there is also the possibility that this leads to a real solution. If the root obtained is real, it is necessary to change the initial condition and expect that this leads to a complex solution, otherwise, it is necessary to change the value of the initial condition again. 

The process described above, it is very similar to what happens when using different  values $ \alpha $ in fractional operators until we find a solution that meets some desired condition. Seeing Newton-Raphson method from the perspective of fractional calculus, one can consider that an order $ \alpha $ remains fixed , in this case $ \alpha = 1 $, and the initial conditions $ x_0 $ are varied until obtaining one solution that satisfies an established criteria. Then reversing the behavior of $ \alpha $ and $ x_0 $, i.e., leave the initial condition fixed and varying the order $ \alpha $, the fractional Newton-Raphson method is obtained, which is nothing other thing than Newton-Raphson method using any definition of fractional derivative that fits the function with which one is working. This change, although in essence simple, allows to find roots in the complex space using real initial conditions,  because fractional operators generally do not carry polynomials to polynomials.

\subsection{Fixed Point Method}

A classic problem of common interest in both Physics and Mathematics is to find the zeros of a function $f:\Omega \subset \nset{R}^n \to \nset{R}^n$, i.e.,

\begin{eqnarray*}
\set{\xi \in \Omega \ : \ f(\xi)=0},
\end{eqnarray*}

this problem often arises as a consequence of wanting to solve other problems, for instance, if we want to determine the eigenvalues of a matrix or want to build a box with a given volume but with minimal surface; in the first example, we need to find the zeros (or roots) of the characteristic polynomial of the matrix, while in the second one we need to find the zeros of the gradient of a function that relates the surface of the box with its volume.

Although finding the zeros of a function may seem like a simple problem, in general, it involves solving nonlinear equations, which in many cases does not have an analytical solution, an example of this is present when we are trying to determine the zeros of the following function

\begin{eqnarray*}
f(x)=\sin(x)-\dfrac{1}{x}.
\end{eqnarray*}

Because in many cases there is no analytical solution, numerical methods are needed to try to determine the solutions to these problems; it should be noted that when using numerical methods, the word \com{determine} should be interpreted as approach a solution with a degree of precision desired. The numerical methods mentioned above are usually of the iterative type and work as follows: suppose we have a function $f : \Omega \subset \nset{R}^n \to \nset{R}^n$ and we search a value $\xi\in \nset{R}^n$ such that $f (\xi) = 0$, then we can start by giving an initial value $x_0\in \nset{R}^n$ and then calculate a value $x_i$ close to the searched value $\xi$ using an iteration function $\Phi : \nset{R}^n \to \nset{R}^n$ as follows \cite{stoer2013}:

\begin{eqnarray}\label{eq:c2.01}
x_{i+1}=\Phi(x_i), & i=0,1,2,\cdots,
\end{eqnarray}

this generates a sequence $\set{x_i}_{i=0}^\infty$, which under certain conditions satisfies that

\begin{eqnarray*}
\lim_{i\to \infty}x_i\to \xi.
\end{eqnarray*}

To understand the convergence of the iteration function $ \Phi $ it is necessary to have the following definition \cite{plato2003concise}:

\begin{definition}
Let $ \Phi: \nset {R} ^ n \to \nset {R} ^ n $ be an iteration function. The method given in \eqref{eq:c2.01} to determine $ \xi \in \nset{R} ^ n $, it is called (locally) \textbf{convergent}, if exists $ \delta> 0 $ such that for all initial value

\begin{eqnarray*}
x_0\in B(\xi;\delta):=\set{y\in \nset{R}^n \ : \ \norm{y-\xi}<\delta},
\end{eqnarray*}

it holds that

\begin{eqnarray}\label{eq:c2.02}
\lim_{i\to \infty}\norm{x_i-\xi}\to 0 & \Rightarrow & \lim_{i\to \infty}x_i=\xi,
\end{eqnarray}

where $ \norm{ \ \cdot \ }: \nset{R} ^ n \to \nset{R} $ denotes any vector norm.

\end{definition}

When it is assumed that the iteration function $ \Phi $ is continuous at $ \xi $ and that the sequence $ \set{x_i} _{i = 0} ^ \infty $ converges to $ \xi $ under the condition given in \eqref{eq:c2.02}, it is true that

\begin{eqnarray}\label{eq:c2.07}
\xi=\lim_{i\to \infty}x_{i+1}=\lim_{i\to \infty}\Phi(x_i)=\Phi\left(\lim_{i\to \infty}x_i\right)=\Phi(\xi),
\end{eqnarray}

the previous result is the reason why the method given in \eqref{eq:c2.01} is called \textbf{Fixed Point Method}.

\subsubsection{Convergence and Order of Convergence}

The (local) convergence of the  iteration function $ \Phi $ established in \eqref{eq:c2.02}, it is useful for demonstrating certain intrinsic properties of the fixed point method. Before continuing it is necessary to have the following definition \cite{stoer2013}:

\begin{definition}
Let $ \Phi: \Omega \subset \nset{R}^n \to \nset{R} ^ n $. The function $ \Phi $ is a \textbf{contraction} on a set $ \Omega_0 \subset \Omega $, if exists a non-negative constant $ \beta <1 $ such that

\begin{eqnarray}\label{eq:c2.03}
\norm{\Phi(x)-\Phi(y)}\leq \beta \norm{x-y}, & \forall  x,y\in \Omega_0,
\end{eqnarray}

where $ \beta $ is called a contraction constant.

\end{definition}

The previous definition guarantee that if the iteration function $ \Phi $ is a contraction on a set $ \Omega_0 $, then it is  Lipschitz continuous and has at least one fixed point. The existence of a fixed point is guaranteed by the following theorem \cite{ortega1990numerical}:

\begin{theorem}
\textbf{Contraction Mapping Theorem}: Let $ \Phi: \Omega \subset \nset{R}^ n \to \nset{R}^n $. Assuming that $ \Phi $ is a contraction on a closed set $ \Omega_0 \subset \Omega $, and that $ \Phi (x) \in \Omega_0 \ \forall x \in \Omega_0 $. Then $ \Phi $ has a unique fixed point $ \xi \in \Omega_0 $ and for any initial value $ x_0 \in \Omega_0 $, the sequence $ \set{x_i}_ {i = 0} ^ \infty $ generated by \eqref{eq:c2.01} converges to $ \xi $. Moreover

\begin{eqnarray}\label{eq:c2.04}
\norm{x_{k+1}-\xi}\leq \dfrac{\beta}{1-\beta} \norm{x_{k+1}-x_{k}}, & k=0,1,2,\cdots,
\end{eqnarray}

where $ \beta $ is the contraction constant given in \eqref{eq:c2.03}.

\end{theorem}

When the fixed point method given by \eqref{eq:c2.01} is used, in addition to convergence, exists a special interest in the \textbf{order of convergence}, which is defined as follows \cite{plato2003concise}:

\begin{definition}
Let $ \Phi: \Omega \subset \nset{R}^ n \to \nset{R}^ n $ be an iteration function with a fixed point $ \xi \in \Omega $. Then the method \eqref{eq:c2.01} is called (locally) \textbf{convergent of (at least) order $ \boldsymbol{p} $} ($ p \geq 1 $), if exists $ \delta> 0 $  and exists a non-negative constant $ C $ (with $ C <1 $ if $ p = 1 $) such that for any initial value $ x_0 \in B (\xi; \delta) $ it is true that

\begin{eqnarray}\label{eq:c2.08}
\norm{x_{k+1}-\xi}\leq C \norm{x_k-\xi}^p, & k=0,1,2,\cdots,
\end{eqnarray}

where $ C $ is called convergence factor.

\end{definition}

The order of convergence is usually related to the speed at which the sequence generated by \eqref{eq:c2.01} converges. For the particular cases $ p = 1 $ or $ p = 2 $ it is said that the method has (at least) linear or quadratic convergence, respectively. The following theorem  for the one-dimensional case \cite{plato2003concise}, allows characterizing the order of convergence of an iteration function $ \Phi $ with its derivatives 

\begin{theorem}\label{teo:c2.01}
Let $ \Phi: \Omega \subset \nset{R} \to \nset {R} $ be an iteration function with a fixed point $ \xi \in \Omega $. Assuming that $\Phi $ is $ p$-times differentiable in $ \xi $ for some $ p \in \nset{N} $, and furthermore

\begin{eqnarray}\label{eq:c2.09}
\left\{
\begin{array}{cc}
 \abs{\Phi^{(k)}(\xi)}=0, \ \forall  k\leq p-1, & \mbox{if }p\geq 2, \\
 \abs{\Phi^{(1)}(\xi)}<1, & \mbox{if }p=1,
\end{array}\right.
\end{eqnarray}

then $ \Phi $ is (locally) convergent of (at least) order $ p $.

\begin{proof}
Assuming that $ \Phi: \Omega \subset \nset{R} \to \nset{R} $ is an iteration function $ p $-times differentiable with a fixed point $ \xi \in \Omega $, then we can expand in Taylor series the function $ \Phi(x_i) $ around $ \xi $ and order $ p $

\begin{eqnarray*}
\Phi(x_i)=\Phi(\xi)+\sum_{s=1}^p \dfrac{\Phi^{(s)}(\xi)}{s!}(x_i-\xi)^s +o \left((x_i-\xi)^p \right),
\end{eqnarray*}

then we obtain 

\begin{eqnarray*}
\abs{\Phi(x_i)-\Phi(\xi)}\leq \sum_{s=1}^p \abs{\dfrac{\Phi^{(s)}(\xi)}{s!}}\abs{x_i-\xi}^s +o \left(\abs{x_i-\xi}^p \right),
\end{eqnarray*}

assuming that the sequence $ \set{x_i}_{i = 0} ^ \infty $ generated by \eqref{eq:c2.01} converges to $ \xi $ and also that $ \abs{\Phi^{(s )}(\xi)} = 0 \ \forall s <p $, the previous expression implies that

\begin{eqnarray*}
\dfrac{\abs{x_{i+1}-\xi}}{\abs{x_i-\xi}^p} =\dfrac{\abs{\Phi(x_{i})-\Phi(\xi)}}{\abs{x_i-\xi}^p}\leq\abs{\dfrac{\Phi^{(p)}(\xi)}{p!}} + \dfrac{o\left(\abs{x_i-\xi}^p \right)}{\abs{x_i-\xi}^p} \underset{i\to \infty}{\longrightarrow} \abs{\dfrac{\Phi^{(p)}(\xi)}{p!}} ,
\end{eqnarray*}

as consequence, exists a value $ k> 0 $ such that

\begin{eqnarray*}
\abs{x_{i+1}-\xi} \leq \abs{\dfrac{\Phi^{(p)}(\xi)}{p!}}\abs{x_i-\xi}^p, & \forall  i\geq k.
\end{eqnarray*}

\end{proof}

\end{theorem}

A version of the previous theorem for the case $n$-dimensional may be found in the reference \cite{stoer2013}.

\subsection{Newton-Raphson Method }

The previous theorem in its n-dimensional form is usually very useful to generate a fixed point method with an order of convergence desired, an order of convergence that is usually appreciated in iterative methods is the (at least) quadratic order. If we have a function $f:\Omega \subset \nset{R}^n \to \nset{R}^n$ and we search a value $\xi\in \Omega$ such that $f(\xi)=0$, we may build an iteration function $ \Phi $ in general form as \cite{burden2002analisis}:

\begin{eqnarray}\label{eq:c2.13}
\Phi(x)=x-A(x)f(x),
\end{eqnarray}

with $ A (x) $ a matrix like

\begin{eqnarray}\label{eq:c2.14}
A(x):=\left( a_{jk}(x)\right)=\begin{pmatrix}
a_{11}(x) & a_{12}(x)& \cdots & a_{1n}(x)\\
a_{12}(x) & a_{22}(x) & \cdots &a_{2n}(x)\\
\vdots & \vdots &\ddots & \vdots \\
a_{n1}(x)& a_{n2}(x)& \vdots & a_{nn}(x)
\end{pmatrix},
\end{eqnarray}

where $a_{jk}:\nset{R}^n \to \nset{R}$  ($1\leq j,k\leq n$). Notice that the matrix $ A (x) $ is determined according to the order of convergence desired. Before continuing, it is necessary to mention that the following conditions are needed:

\begin{enumerate}
\item Suppose we can generalize the \mref{Theorem}{\ref{teo:c2.01}} to the case $ n $-dimensional, although for this it is necessary to guarantee that the iteration function $ \Phi $ given by \eqref{eq:c2.13} near the value $ \xi $ can be expressed  in terms of its Taylor series in several variables.
\item It is necessary to guarantee that the norm of the equivalent of the first derivative in several variables of the iteration function $ \Phi $ tends to zero near the value $ \xi $. 
\end{enumerate}

Then, we will assume that the first condition is satisfied; for the second condition we have that the equivalent to the first derivative in several variables is the \textbf{Jacobian matrix} of the function $\Phi$, which is defined as follows \cite{ortega1990numerical}:

\begin{eqnarray}\label{eq:c2.15}
\Phi^{(1)}(x):=\left(\partial_j \Phi_k(x) \right) ,
\end{eqnarray}

where

\begin{eqnarray*}
\partial_j\Phi_k(x):= \dfrac{\partial }{\partial x_j}\Phi_k(x), &1\leq j,k\leq n,
\end{eqnarray*}

with $\Phi_k:\nset{R}^n \to \nset{R}$, the competent $ k$-th of the iteration function $ \Phi $. Now considering that 

\begin{eqnarray}\label{eq:c2.16}
\begin{array}{cccc}
\ds \lim_{x\to \xi} \norm{\Phi^{(1)}(x)}=0 & \Rightarrow &\ds  \lim_{x\to \xi}\partial_j \Phi_k(x)=0, & \forall j,k\leq n,
\end{array}
\end{eqnarray}

Assume that we have a function $f(x):\Omega \subset \nset{R}^n \to \nset{R}^n$  with a zero $\xi \in \Omega$, such that all of its first partial derivatives are defined in $ \xi $. Then taking the iteration function $ \Phi $ given by \eqref{eq:c2.13}, the $ k $-th component of the iteration function may be written as

\begin{eqnarray*}
\Phi_k(x)=x_k-\sum_{j=1}^na_{kj}(x)f_j(x),
\end{eqnarray*}

then

\begin{eqnarray*}
\partial_l \Phi_k(x)=\delta_{lk}-\sum_{j=1}^n \left[ a_{kj}(x)\partial_l f_j(x)+\left(\partial_l a_{kj}(x) \right)f_k(x) \right],
\end{eqnarray*}

where $ \delta_{lk} $ is the Kronecker delta, which is defined as

\begin{eqnarray*}
\delta_{lk}=\left\{
\begin{array}{cc}
1,& \mbox{si }l=k,\\
0,& \mbox{si }l\neq k.
\end{array}\right.
\end{eqnarray*}

Assuming that \eqref{eq:c2.16} is fulfilled

\begin{eqnarray*}
\partial_l \Phi_k(\xi)=\delta_{lk}-\sum_{j=1}^n a_{kj}(\xi)\partial_l f_j(\xi)=0 & \Rightarrow & \sum_{j=1}^n a_{kj}(\xi)\partial_l f_j(\xi)=\delta_{lk}, \ \forall l,k\leq n,
\end{eqnarray*}

the previous expression may be written in matrix form as

\begin{eqnarray*}
A(\xi) f^{(1)}(\xi)=I_n & \Rightarrow & A(\xi)= \left(f^{(1)}(\xi)\right)^{-1},
\end{eqnarray*}

where $ f ^{(1)} $ and $ I_n $ are the Jacobian matrix of the  function $ f $ and the identity matrix of $ n \times n $, respectively. Denoting by $ \dt{A} $ the determinant of the  matrix $ A $, then any matrix $ A (x) $ that fulfill the following condition

\begin{eqnarray}\label{eq:c2.17}
\lim_{x\to \xi}A(x)= \left(f^{(1)}(\xi)\right)^{-1}, & \dt{f^{(1)}(\xi)}\neq 0,
\end{eqnarray}

guarantees that exists $ \delta> 0 $ such that the iteration function $ \Phi $ given by \eqref{eq:c2.13}, converges (locally) with an order of convergence (at least) quadratic in $ B (\xi; \delta) $. The following fixed point method can be obtained from the previous result

\begin{eqnarray}\label{eq:c2.18}
x_{i+1}:=\Phi(x_i)=x_i-\left(f^{(1)}(x_i) \right)^{-1}f(x_i) , & i=0,1,2,\cdots,
\end{eqnarray}

which is known as  Newton-Raphson method ($n$-dimensional), also known as Newton's method \cite{ortega1970iterative}.

Although the condition given in \eqref{eq:c2.17} could seem that Newton-Raphson method always has an order of convergence (at least) quadratic, unfortunately, this is not true; the order of convergence is now conditioned by the way in which the function f is constituted, the mentioned above may be appreciated in the following proposition

\begin{proposition}
Let $ f: \Omega \subset \nset {R} \to \nset{R} $ be a function that is at least twice differentiable in $ \xi \in \Omega $. So if $ \xi $ is a zero of $ f $ with algebraic multiplicity $ m $ ($ m \geq 2 $), i.e.,

\begin{eqnarray*}
f(x)=(x-\xi)^m g(x), & g(\xi)\neq 0,
\end{eqnarray*}

the Newton-Raphson method (one-dimensional) has an order of convergence (at least) linear.

\begin{proof}
Suppose we have $ f: \Omega \subset \nset{R} \to \nset{R} $ a function with a zero $ \xi \in \Omega $ of algebraic multiplicity $ m \geq 2 $, and that $ f $ is at least twice differentiable in $ \xi $, then

\begin{align*}
f(x)=&(x-\xi)^mg(x), \ \  \ g(\xi)\neq0, \\
f^{(1)}(x)=& (x-\xi)^{m-1}\left(m g(x)+(x-\xi)g^{(1)}(x)\right),
\end{align*}

as consequence the derivative of the iteration function $ \Phi $ of Newton-Raphson method may be expressed as

\begin{eqnarray*}
\begin{array}{c}
\Phi^{(1)}(x)=1-\dfrac{mg^2(x)+(x-\xi)^2\left[\left(g^{(1)} (x)\right)^2-g(x)g^{(2)}(x) \right]}{\left(mg(x)+(x-\xi)g^{(1)}(x) \right)^2},
\end{array}
\end{eqnarray*}

therefore

\begin{eqnarray*}
\lim_{x\to \xi}\Phi^{(1)}(x)= 1-\dfrac{1}{m},
\end{eqnarray*}

and by the \mref{Theorem}{\ref{teo:c2.01}}, then Newton-Raphson method  under the hypothesis of the proposition converges (locally) with an order of convergence (at least) linear.

\end{proof}

\end{proposition}

Notice that the order of convergence (at least) quadratic, does not necessarily imply a higher speed of convergence when using Newton-Raphson method; to verify this, just is necessary to take a polynomial of the form

\begin{eqnarray*}
f(x)=x^n+\sum_{k=0}^{n-1} a_kx^k, & a_k\in \nset{R}, \ \forall k\leq n-1,
\end{eqnarray*}

now assuming that $ \xi \in \nset{R} $ is a zero of the polynomial and that the sequence $ \set {x_i} _ {i = 0} ^ \infty $ generated by \eqref{eq:c2.18} has an initial condition $ x_0 \in \nset{R} $ away from $ \xi $, then

\begin{eqnarray*}
\ds x_{i+1}=\Phi(x_i)=x_i-\dfrac{ x_i^n+\sum_{k=0}^{n-1} a_kx_i^k}{ nx_i^{n-1}+\sum_{k=1}^{n-1} ka_kx_i^{k-1}}\approx x_i\left(1-\dfrac{1}{n} \right),
\end{eqnarray*}

so a small change is obtained between the values $ x_i $ and $ x_{i + 1} $; therefore, the sequence generated by Newton method will have a slow convergence to the value of $ \xi $ \cite{stoer2013}.

\section{Basic Definitions of the Fractional Derivative}

\subsection{Introduction to the Definition of Riemann-Liouville }

One of the key pieces in the study of fractional calculus is the iterated integral, which is defined as follows \cite{hilfer00}:

\begin{definition}
Let $ L_{loc} ^ 1 (a, b) $, the space of locally integrable functions in the interval $ (a, b) $. If $ f $ is a function such that $ f \in L_ {loc} ^ 1 (a, \infty) $, then the $n$-th iterated integral of the function $ f $ is given by \cite{hilfer00}:

\begin{eqnarray}\label{eq:c1.16}
\begin{array}{c}
\ds \ifr{}{a}{I}{x}{n} f(x)=\ifr{}{a}{I}{x}{}\left(\ifr{}{a}{I}{x}{n-1} f(x)  \right)=\frac{1}{(n-1)!}\int_a^x(x-t)^{n-1}f(t)dt,
\end{array}
\end{eqnarray}

where

\begin{eqnarray*}
\ifr{}{a}{I}{x}{} f(x):=\int_a^x f(t)dt.
\end{eqnarray*}

\end{definition}

Considerate that $ (n-1)! = \gam{n} $
, a generalization of \eqref{eq:c1.16} may be obtained for an arbitrary order $ \alpha> 0 $

\begin{eqnarray}\label{eq:c1.17}
\ifr{}{a}{I}{x}{\alpha} f(x)=\dfrac{1}{\gam{\alpha}}\int_a^x(x-t)^{\alpha-1}f(t)dt,
\end{eqnarray}

similarly, if $ f \in L_{loc} ^ 1 (- \infty, b) $, we may define

\begin{eqnarray}\label{eq:c1.18}
\ifr{}{x}{I}{b}{\alpha} f(x)=\dfrac{1}{\gam{\alpha}}\int_x^b(t-x)^{\alpha -1}f(t)dt,
\end{eqnarray} 

the equations \eqref{eq:c1.17} and \eqref{eq:c1.18} correspond to the definitions of \textbf{right and left fractional integral of Riemann-Liouville}, respectively.

Below, there are some valid properties for right fractional integrals. However, the same properties for left fractional integrals have some changes. Fractional integrals satisfy the  \textbf{semigroup property}, which is given in the following proposition \cite{hilfer00}:

\begin{proposition}
Let $ f $ be a function. If $ f \in L_{loc} ^ 1 (a, \infty) $, then the fractional integrals of $ f $ satisfy that

\begin{eqnarray}\label{eq:c1.19}
\ifr{}{a}{I}{x}{\alpha} \ifr{}{a}{I}{x}{\beta}f(x) = \ifr{}{a}{I}{x}{\alpha + \beta}f(x),& \alpha,\beta>0.
\end{eqnarray}

\end{proposition}

From the previous result, it is obtained that in particular

\begin{eqnarray*}
\ifr{}{a}{I}{x}{n}\ifr{}{a}{I}{x}{\alpha}f(x)=\ifr{}{a}{I}{x}{n+\alpha}f(x),& n\in\mathds{N}, \ \alpha>0,
\end{eqnarray*}

since the $ d / dx $ operator is the inverse operator
to the left of the operator $ \ifr {}{a}{I}{x}{} $, any integral $ \alpha$-th of a function $ f \in L_{loc} ^ 1 (a, \infty) $ may be written as

\begin{eqnarray}\label{eq:c1.20}
\ifr{}{a}{I}{x}{\alpha}f(x)=\dfrac{d^n}{dx^n}\left( \ifr{}{a}{I}{x}{n}\ifr{}{a}{I}{x}{\alpha}f(x) \right)=\dfrac{d^n}{dx^n}\left( \ifr{}{a}{I}{x}{n+\alpha}f(x)\right).
\end{eqnarray}

It should be mentioned that the above results are also valid for $\alpha$ in complexes, with $ \re{\alpha}> 0 $. Under appropriate conditions in \eqref{eq:c1.20} its follows 

\begin{eqnarray}\label{eq:c1.21}
\ifr{}{a}{I}{x}{-n}f(x)=\dfrac{d^n}{dx^n}\left( \ifr{}{a}{I}{x}{0}f(x)\right)=\der{d}{x}{n}f(x),
\end{eqnarray}

then, from \eqref{eq:c1.20} and \eqref{eq:c1.21} the following operator is built, which corresponds to the \textbf{ (right) fractional derivative  of Riemann-Liouville}

\begin{eqnarray}
\ifr{}{a}{D}{x}{\alpha}f(x):=\ifr{}{a}{I}{x}{-\alpha}f(x)=\dfrac{d^n}{dx^n}\left( \ifr{}{a}{I}{x}{n-\alpha}f(x)\right)=\frac{1}{\gam{n-\alpha} }\der{d}{x}{n}\int_a^x (x-t)^{n-\alpha -1} f(t)dt, \label{eq:c1.22}
\end{eqnarray}

where $ n = \lfloor \re{\alpha} \rfloor + 1 $, with $ \lfloor \re{\alpha} \rfloor $ the largest integer less than or equal to $ \re{\alpha} $.  We can unify the definitions of integral and fractional derivative of Riemann-Liouville, given by \eqref{eq:c1.17} and \eqref{eq:c1.22}, as follows

\begin{eqnarray}\label{eq:c1.23}
\begin{array}{c}
\ifr{}{a}{D}{x}{\alpha}f(x) := \left\{
\begin{array}{cc}
\ds \ifr{}{a}{I}{x}{-\alpha}f(x), &\mbox{if }\re{\alpha}<0,\\  
\ds \dfrac{d^n}{dx^n}\left( \ifr{}{a}{I}{x}{n-\alpha}f(x)\right), & \mbox{if }\re{\alpha}\geq 0, 
\end{array}
\right.
\end{array}
\end{eqnarray}

where  $ n = \lfloor \re{\alpha} \rfloor + 1 $. Applying the  operator \eqref{eq:c1.23} with $ a = 0 $ and $ \alpha \in \nset{R} \setminus \nset {Z} $ to the  function $ x^{\mu} $, we obtain

\begin{eqnarray}\label{eq:c1.13}
\ifr{}{0}{D}{x}{\alpha}x^\mu = \left\{
\begin{matrix}
(-1)^\alpha \dfrac{\gam{-\left( \mu + \alpha \right) }}{\gam{-\mu}}x^{\mu-\alpha}, & \mbox{if }\mu\leq -1,\\
 \dfrac{\gam{\mu+1}}{\gam{\mu-\alpha+1}}x^{\mu-\alpha}, & \mbox{if } \mu>-1.
\end{matrix}
 \right.
\end{eqnarray}

\subsection{Introduction to the Definition of Caputo }

Michele Caputo (1969) published a book and introduced a new definition of fractional derivative, he created this definition with the objective of modeling anomalous diffusion phenomena. The definition of Caputo had already been discovered independently by Gerasimov (1948). This fractional derivative is of the utmost importance since it allows us to give a physical interpretation of the initial value problems, moreover to being used to model fractional time. In some texts, it is known as the fractional derivative of Gerasimov-Caputo \cite{kilbas2006theory}.

Let $ f $ be a function, such that $ f $ is $ n$-times differentiable with $ f ^{(n)} \in L_{loc}^ 1 (a, b) $, then the \textbf{(right) fractional derivative  of Caputo} is defined as:

\begin{align}
\ifr{C}{a}{D}{x}{\alpha}f(x):= &\ifr{}{a}{I}{x}{n-\alpha}\left( \der{d}{x}{n} f(x)\right) = \dfrac{1}{\gam{n-\alpha}}\int_{a}^{x} (x-t)^{n-\alpha -1} f^{(n)}(t)dt , \label{eq:c1.25}
\end{align}

where $n=\lfloor \re{\alpha}\rfloor+1$. It should be mentioned that the fractional derivative of Caputo behaves as the inverse operator to the left of fractional integral of Riemann-Liouville , that is,

\begin{eqnarray*}
\ifr{C}{a}{D}{x}{\alpha}(\ifr{}{a}{I}{x}{\alpha}f(x))=f(x),
\end{eqnarray*}

On the other hand, the relation between the fractional derivatives of Caputo and Riemann-Liouville is given by the following expression \cite{kilbas2006theory}

\begin{eqnarray*}
\ifr{C}{a}{D}{x}{\alpha}f(x)=\ifr{}{a}{D}{x}{\alpha}\left(f(x)-\sum_{k=0}^{n-1}\dfrac{f^{(k)}(a)}{k!}(x-a)^k\right), 
\end{eqnarray*}

then if $f^{(k)}(a)=0 \ \ \forall k<n$, we obtain

\begin{eqnarray*}
\ifr{C}{a}{D}{x}{\alpha}f(x)=\ifr{}{a}{D}{x}{\alpha}f(x), 
\end{eqnarray*}

considering the previous particular case, it is possible to unify the definitions of fractional integral of Riemann-Liouville and  fractional derivative of Caputo as follows

\begin{eqnarray}\label{eq:c1.233}
\begin{array}{c}
\ifr{C}{a}{D}{x}{\alpha}f(x) := \left\{
\begin{array}{cc}
\ds \ifr{}{a}{I}{x}{-\alpha}f(x), &\mbox{si }\re{\alpha}<0,\\  
\ds \ifr{}{a}{I}{x}{n-\alpha}\left( \der{d}{x}{n} f(x)\right) , & \mbox{si }\re{\alpha}\geq 0, 
\end{array}
\right.
\end{array}
\end{eqnarray}

\section{Fractional Newton Method}

Let $\nset{P}_n(\nset{R})$, the space of polynomials of degree less than or equal to $ n $ with real coefficients. The zeros $ \xi $ of a function $ f \in \nset {P} _n (\nset{R}) $ are usually named as roots. The Newton-Raphson method is useful for finding the roots of a function $ f$. However, this method is limited because it cannot find roots $ \xi \in \nset{C} \setminus \nset {R} $, if the sequence $ \set{x_i}_{i = 0} ^ \infty $ generated by \eqref{eq:c2.18} has an initial condition $ x_0 \in \nset{R} $. To solve this problem and develop a method that has the ability to find roots, both real and complex, of a polynomial if the initial condition $ x_0 $ is real, we propose a new method called fractional Newton-Raphson method, which consists of Newton-Raphson method with the implementation of the fractional derivative. Before continuing, it is necessary to define the \textbf{fractional Jacobian matrix} of a function $ f: \Omega \subset \nset{R} ^ n \to \nset {R}^ n $ as follows

\begin{eqnarray}\label{eq:c2.26}
f^{(\alpha)}(x):=\left(\partial_j^\alpha f_k(x) \right),
\end{eqnarray}

where

\begin{eqnarray*}
\partial_j^\alpha f_k(x):= \der{\partial}{x_j}{\alpha}f_k(x), &1\leq j,k\leq n.
\end{eqnarray*}

with $f_k:\nset{R}^n \to \nset{R}$. The  operator $ \partial^ \alpha / \partial x_j ^ \alpha $ denotes any fractional derivative, applied only to the variable $ x_j $, that satisfies the following condition of continuity respect to the order of the derivative

\begin{eqnarray*}
\lim_{\alpha \to 1}\der{\partial}{x_j}{\alpha}f_k(x)=\dfrac{\partial}{\partial x_j}f_k(x), & 1\leq j,k\leq n,
\end{eqnarray*}

then, the matrix \eqref{eq:c2.26} satisfies that

\begin{eqnarray}\label{eq:c2.27}
\lim_{\alpha\to 1}f^{(\alpha)}(x)=f^{(1)}(x),
\end{eqnarray}

where $ f ^{(1)} (x) $ denotes the Jacobian matrix of the  function $ f $.

Taking into account that a polynomial of degree $ n $ it is composed of $ n + 1 $ monomials of the form $ x ^ m $, with $m\geq 0$, we can take the equation \eqref{eq:c1.13} with \eqref{eq:c2.18}, to define the following iteration function that results in the \textbf{Fractional Newton-Raphson Method} \cite{fernando2017fractional,brambila2018fractional}

\begin{eqnarray}\label{eq:c2.200}
\begin{array}{cc}
x_{i+1}:= \Phi\left( \alpha, x_i \right)=x_i-\left(f^{\left(\alpha\right)}(x_i) \right)^{-1}  f(x_i),& i=0,1,2,\cdots.
\end{array}
\end{eqnarray}

To try to guarantee that the sequence $ \set{x_i}_ {i = 0} ^\infty $ generated by \eqref {eq:c2.200} has a order of convergence (at least) quadratic, the condition \eqref{eq:c2.17} is combined with \eqref{eq:c2.27} to define the following function

\begin{eqnarray}\label{eq:c2.21}
\alpha_{f}(x):=\left\{
\begin{array}{cc}
\alpha, & \mbox{if } \norm{f(x)}\geq \delta \mbox{ and } \norm{x}\neq 0 ,\\
1, & \mbox{if } \norm{f(x)}<\delta \mbox{ or }\norm{x}=0,
\end{array}\right.
\end{eqnarray}

then,  in \eqref{eq:c2.200} we make $\alpha \to \alpha_f (x) $ obtaining the following iteration function that results in the \textbf{Fractional Newton Method}

\begin{eqnarray}\label{eq:c2.20}
\begin{array}{cc}
x_{i+1}:= \Phi\left( \alpha_{f}(x_i), x_i \right)=x_i-\left(f^{\left(\alpha_{f}(x_i) \right)}(x_i) \right)^{-1}  f(x_i),& i=0,1,2,\cdots,
\end{array}
\end{eqnarray}

the difference between the methods \eqref{eq:c2.200} and \eqref{eq:c2.20},  is that the just for the second  can exists $\delta> 0 $ such that if the sequence $ \set{x_i}_{i = 0} ^ \infty $ generated by \eqref{eq:c2.20} converges to a root $ \xi $ of $ f $,  exists $ k> 0 $ such that $ \forall i \geq k $, the sequence has an order of convergence (at least) quadratic in $B(\xi;\delta)$. 

The value of $ \alpha $ in \eqref{eq:c2.21} is assigned with the following reasoning: when we use the definitions of fractional derivatives given by \eqref{eq:c1.23} and \eqref{eq:c1.233} in a function $ f $, it is necessary that the function be $ n$-times integrable and $n$-times differentiable, where $ n = \lfloor \re{\alpha} \rfloor + 1 $, therefore $ \abs{\alpha} <n $ and, on the other hand, for using Newton method it is just necessary that the function be one-time differentiable, as a consequence of $ \eqref{eq:c2.21} $ it is obtained that

\begin{eqnarray}\label{eq:c2.22}
-2<\alpha<2, & \alpha \neq -1,0,1.
\end{eqnarray}

Without loss of generality, to understand why the sequence $ \set{x_i}_{i = 0}^\infty $ generated by the fractional Newton method when we use a function $ f \in P_ {n}  (\nset {R }) $,   has the ability to enter the complex space starting from an initial condition $ x_0 \in \nset {R} $, it is only necessary to observe the fractional derivative of  Riemann-Liouville of  order $ \alpha = 1/2 $ of the monomial $ x ^ m$

\begin{eqnarray*}
\ifr{}{0}{D}{x}{\frac{1}{2}}x^m = \frac{\sqrt{\pi}}{2\gam{m+\frac{1}{2}}}x^{m-\frac{1}{2}}, & m\geq0,
\end{eqnarray*}

whose result is a function with a rational exponent, contrary to what would happen when using the conventional derivative. When the iteration function given by \eqref{eq:c2.20} is used,  we must taken an initial condition $ x_0 \neq 0 $, as a consequence of the fact that the fractional derivative of order $\alpha > 0 $ of a constant is usually proportional to the function $ x ^{-\alpha}$.

When we take $ \alpha \neq 1 $, the sequence $ \set{x_i} _{i = 0} ^ \infty $ generated by the iteration function \eqref{eq:c2.20}, presents among its behaviors, the following particular cases depending on the initial condition $ x_0 $: \label{pg:c2.01}

\begin{enumerate}
\item If we take an initial condition $ x_0> 0 $, the sequence $ \set{x_i}_{i = 0} ^ \infty $ may be divided into three parts, this occurs because it may exists a value $ M \in \nset{N} $  for which $ \set{x_i}_{i = 0}^{M-1} \subset \nset{R}^+ $ with $ \set{x_M} \subset \nset{R}^{-}$, in consequence $ \set {x_i} _ {i \geq M + 1} \subset \nset{C} $.

\item On the other hand, if  we take an initial  $ x_0 <0 $ condition,  the sequence $ \set {x_i}_ {i = 0} ^ \infty $ may  be divided into two parts, $\set{x_0}\subset \nset{R}^{-}$ and  $\set{x_i}_{i\geq 1}\subset \nset{C}$.
\end{enumerate}

Unlike classical Newton-Raphson method; which uses tangent lines to generate a sequence $ \set{x_i}_{i = 0}^\infty $, the  fractional Newton method uses lines more similar to secants (see Figure \ref{fig:01}). A consequence of the fact that the lines are not tangent when using \eqref{eq:c2.20}, is that different trajectories can be obtained for the same initial condition $x_0 $ just by changing the order $ \alpha $ of the derivative (see Figure \ref{fig:02}).

\begin{figure}[!ht]
\centering
\includegraphics[width=0.6\textwidth, height=0.3\textwidth]{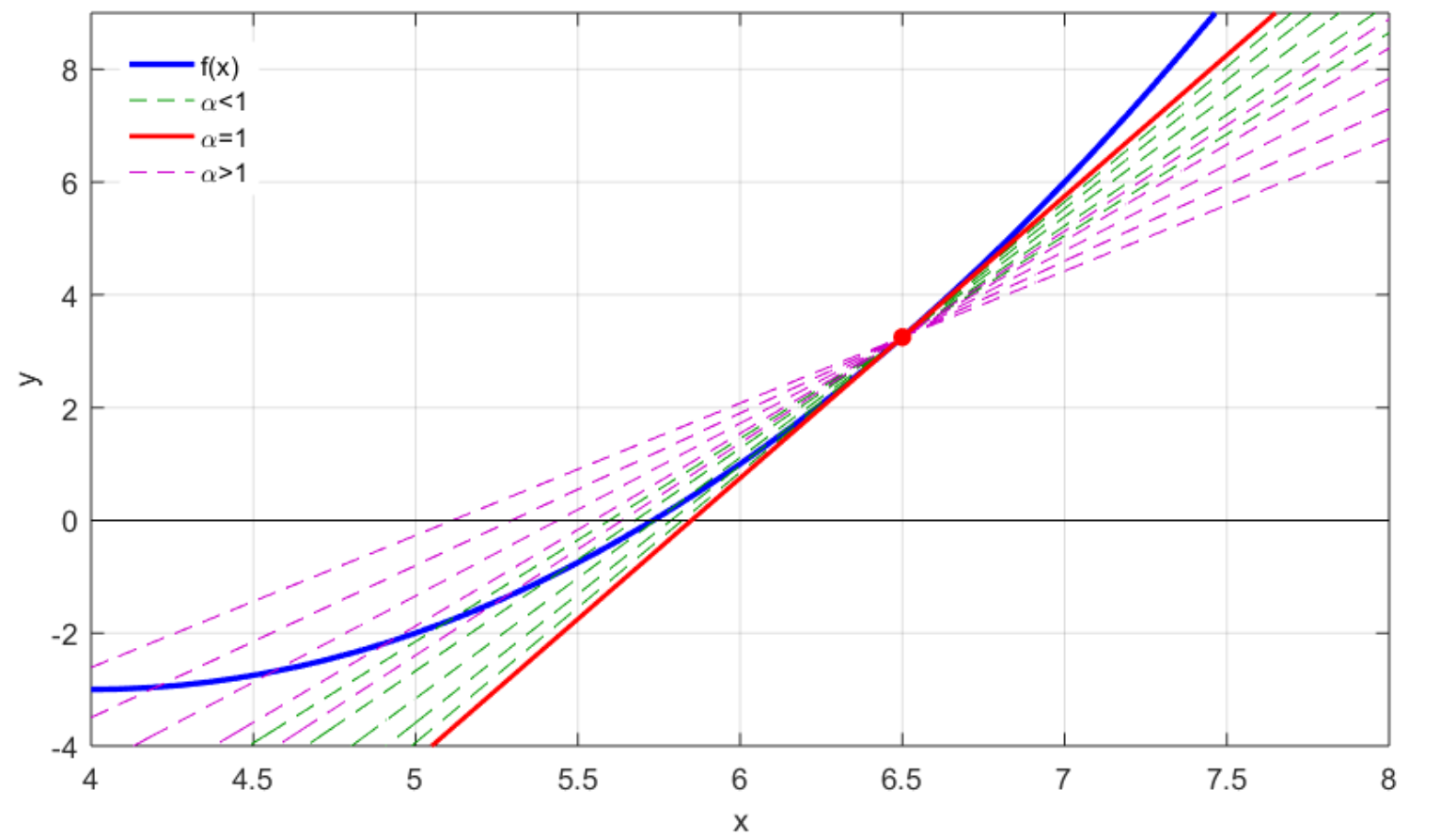}
\caption{ Illustration of some lines generated by the fractional Newton method, the red line corresponds to the case $\alpha = 1 $.  }\label{fig:01}
\end{figure}

\begin{figure}[!ht]
        \begin{subfigure}[c]{0.24\textwidth}
        \centering
 \includegraphics[width=\textwidth, height=0.75\textwidth]{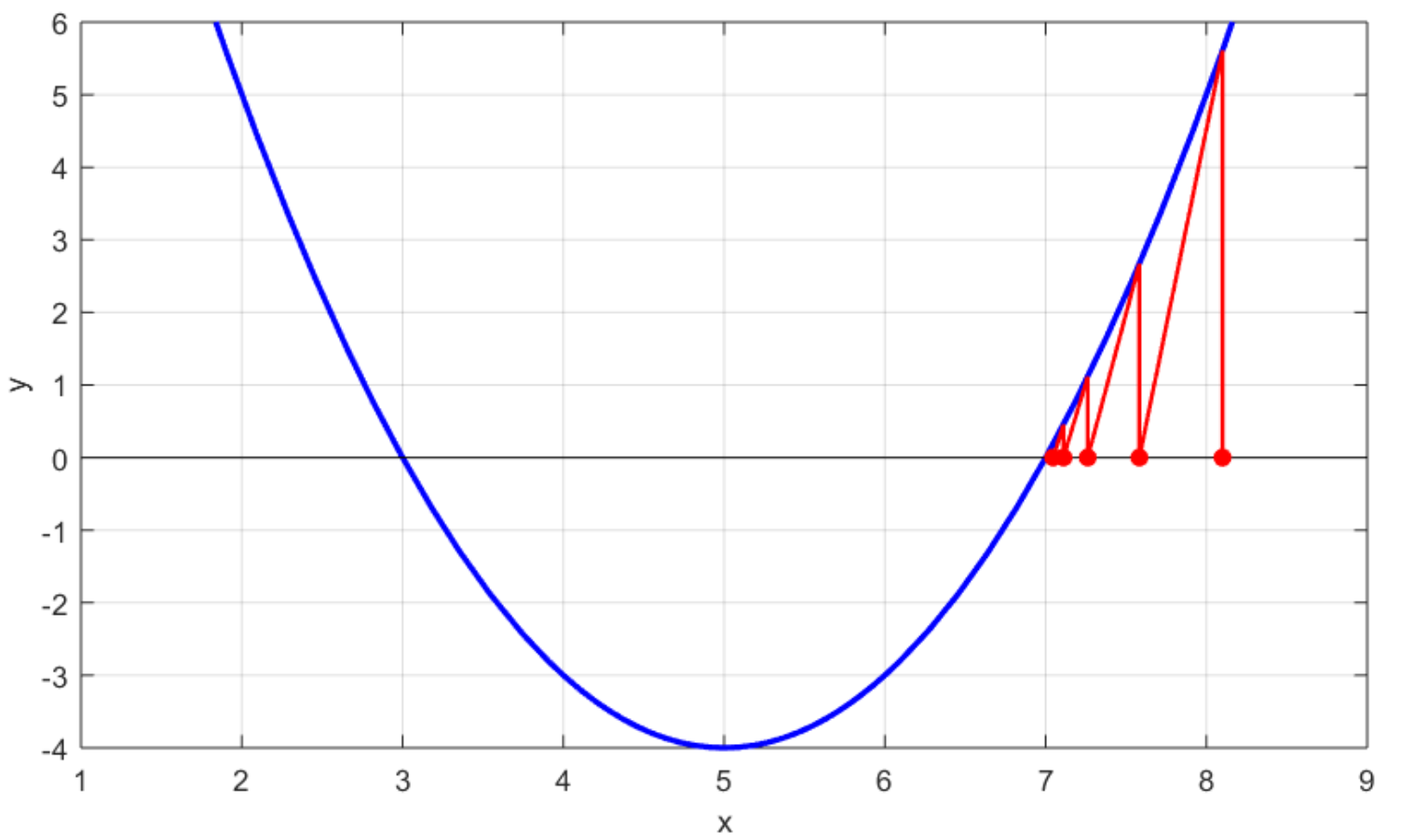}      
    \caption*{a) $\alpha=-0.77$}
    \end{subfigure}
        \begin{subfigure}[c]{0.24\textwidth}
        \centering
 \includegraphics[width=\textwidth, height=0.75\textwidth]{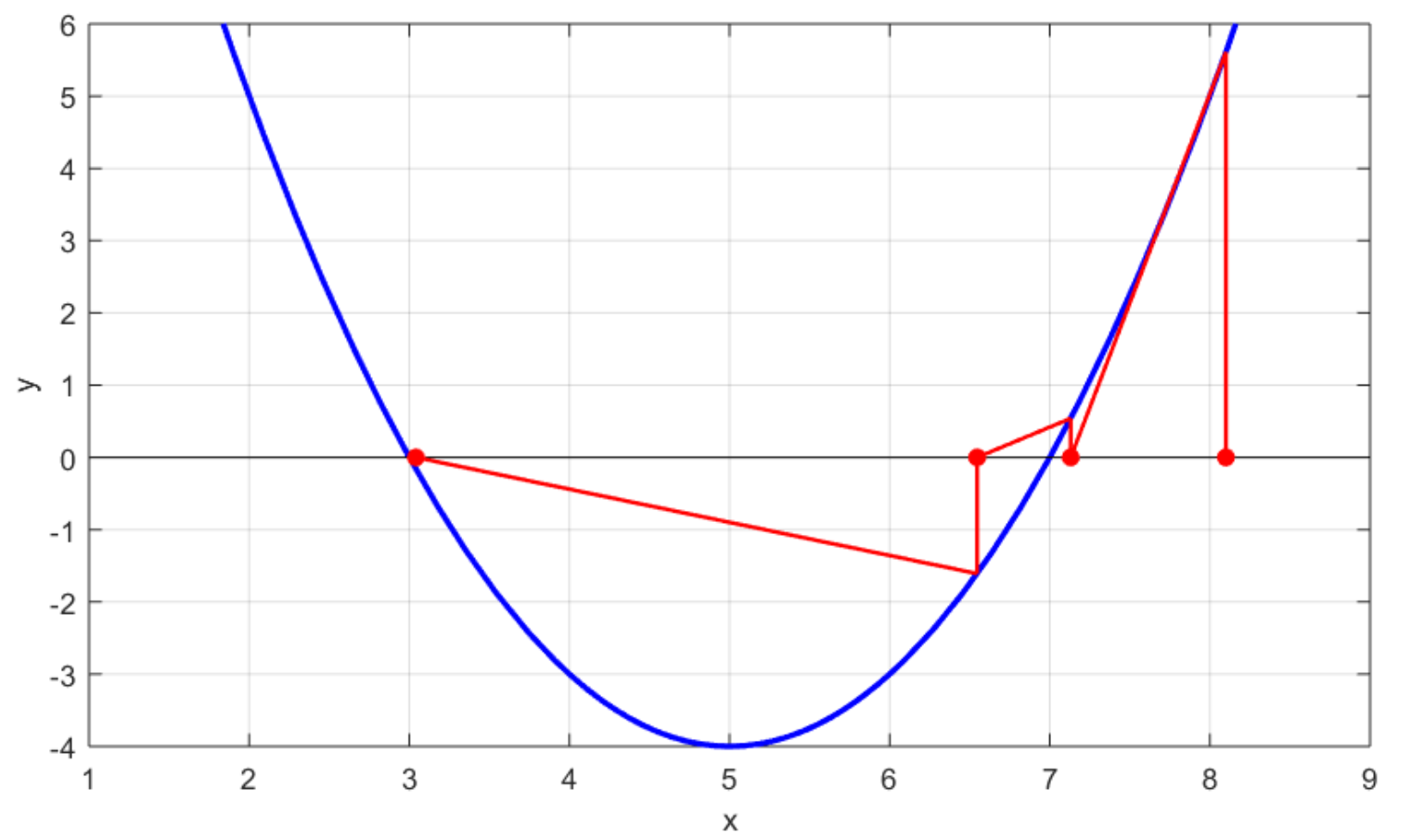}     
    \caption*{b) $\alpha=-0.32$}
    \end{subfigure}    
    \centering
    \begin{subfigure}[c]{0.24\textwidth}
    \centering
 \includegraphics[width=\textwidth, height=0.75\textwidth]{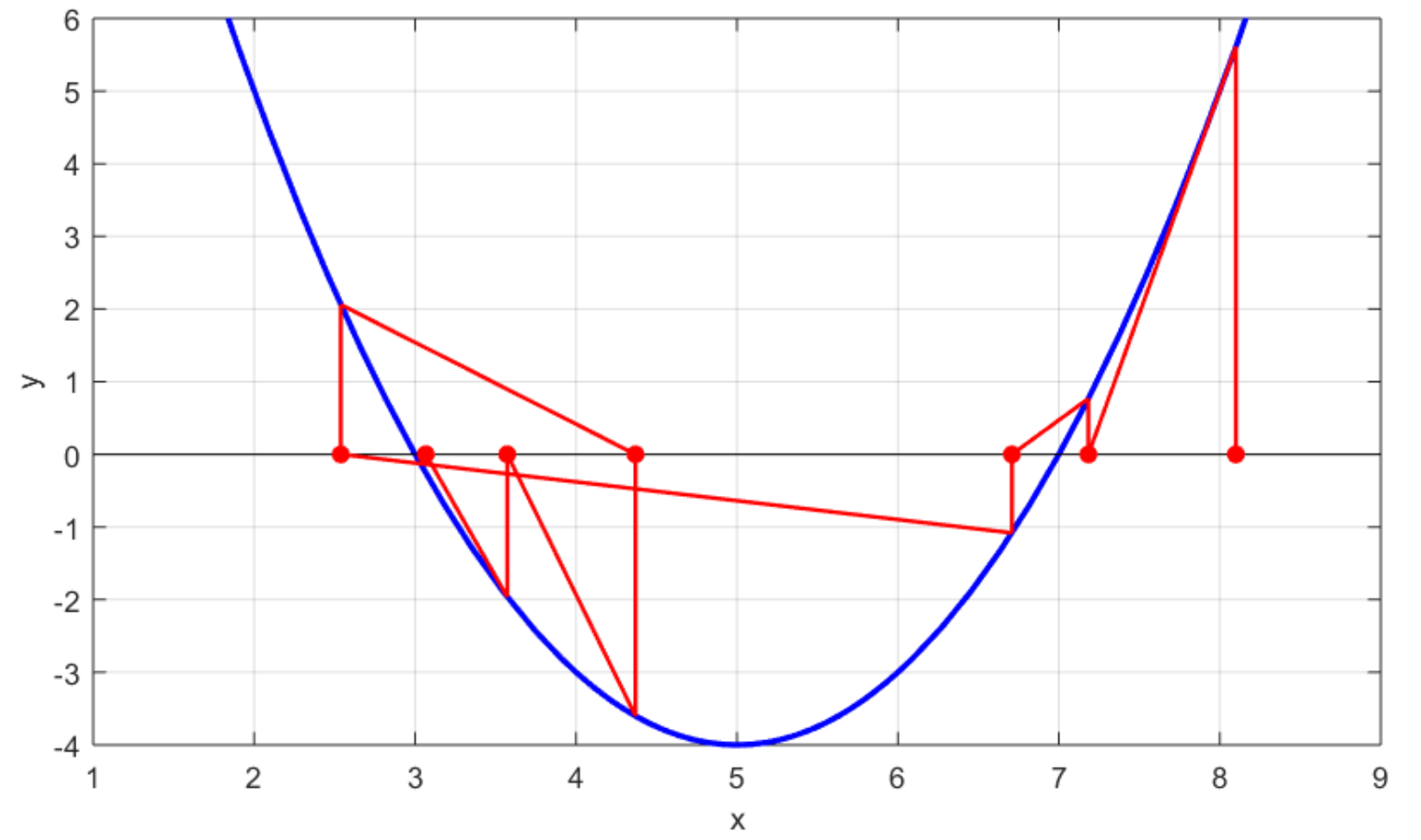}      
    \caption*{c) $\alpha=0.19$}
    \end{subfigure}
    \begin{subfigure}[c]{0.24\textwidth}
    \centering
 \includegraphics[width=\textwidth, height=0.75\textwidth]{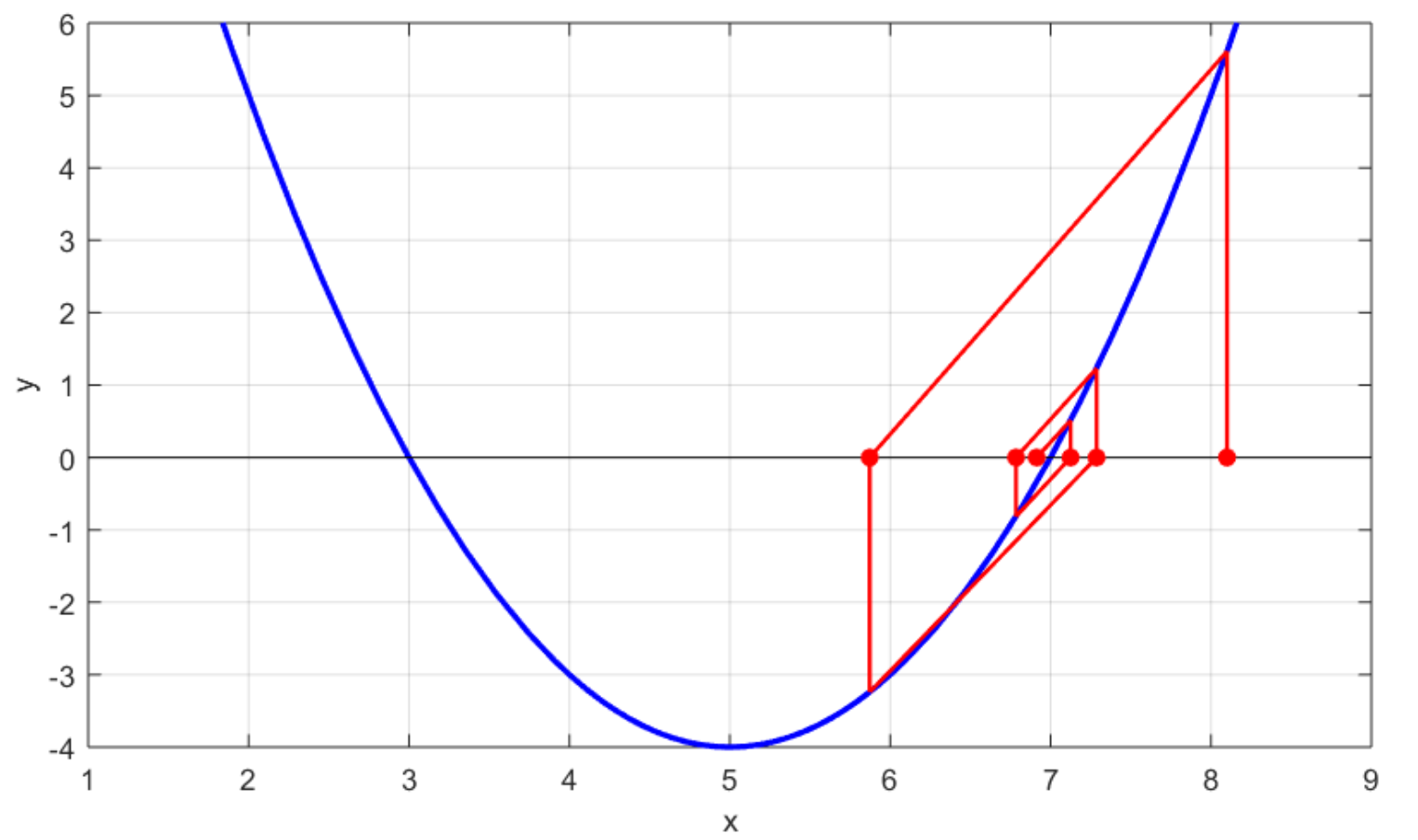}      
    \caption*{d) $\alpha=1.87$}
    \end{subfigure}        
        \caption{llustrations of some trajectories generated by the fractional Newton method for the same initial condition $ x_0 $ but with different values of $ \alpha $.}\label{fig:02}
\end{figure}

A typical inconvenience that arises in problems related to fractional calculus, it is the fact that it is not always known what is the appropriate order $\alpha$   to solve these problems. As a consequence, different values of $\alpha$ are generally tested and we choose  the value that allows to find the best solution considering an criterion of precision established. Based on the aforementioned, it is necessary  to follow the instructions below when using the fractional Newton method to find the zeros $ \xi $ of a  function $ f $: \label{pg:c2.02}

\begin{enumerate}
\item Without considering the integers $ -1, \ 0 $ and $ 1 $, a partition of the interval $ [- 2,2] $ is created as follows

\begin{eqnarray*}
-2=\alpha_{0}<\alpha_1<\alpha_2<\cdots < \alpha_{s}<\alpha_{s+1}=2,
\end{eqnarray*}

and using the partition, the following sequence $ \set{\alpha_m}_ {m = 1} ^ s $ is created.

\item We choose a non-negative tolerance $ TOL <1 $, a limit of iterations $ L_{IT}> 1 $ for all $ \alpha_m $, an initial condition $ x_0 \neq 0 $ and a value $ M>L_{IT} $.

\item We choose a value $ \delta> 0 $ to use $ \alpha_f (x) $ given by \eqref{eq:c2.21}, such that $ TOL <\delta <1 $. In addition, it is taken a fractional derivative  that satisfies the condition of continuity \eqref{eq:c2.27}, and it is unified with the fractional integral in the same way as in the equations \eqref{eq:c1.23} and \eqref{eq:c1.233}.

\item The iteration function \eqref{eq:c2.20} is used with all the values of the partition $ \set {\alpha_m} _ {m = 1} ^ s $, and for each value $ \alpha_m $ is generated a sequence $ \set {{} ^ mx_i} _ {i = 0}^{R_m} $, where

\begin{eqnarray*}
\begin{array}{c}
R_m=\left\{
\begin{array}{cl}
 K_1\leq L_{IT},&  \mbox{if exists $k>0$ such that } \norm{f\left({}^mx_k \right)}\geq M \ \ \forall k\geq i,\\
K_2\leq L_{IT},&
 \mbox{if exists $ k> 0 $ such that } \norm{f\left({}^mx_k \right)}\leq TOL \ \ \forall k\geq i.\\
L_{IT}, & \mbox{if } \norm{f\left({}^mx_i \right)}> TOL \ \ \forall i\geq 0, 
\end{array}\right.
\end{array}
\end{eqnarray*}

then is generated a sequence $ \set{x_ {R_ {m_k}}} _ {k = 1} ^ r $ , with $ r \leq s $, such that

\begin{eqnarray*}
\norm{f\left(x_{R_{m_k}} \right)}\leq TOL , & \forall k\geq 1.
\end{eqnarray*}

\item We choose a value $ \varepsilon> 0 $  and we take the values $x_{R_{m_1}} $ and $ x_ {R_{m_2}} $, then is defined $ X_1 = x_ {R_ {m_1}} $. If the following condition is fulfilled

\begin{eqnarray}\label{eq:c2.23}
\dfrac{\norm{X_1-x_{R_{m_2}}}}{\norm{X_1}} \leq \varepsilon & \mbox{and} & R_{m_2}\leq R_{m_1},
\end{eqnarray}

is taken $ X_1 = x_ {R_ {m_2}} $. On the other hand if

\begin{eqnarray}\label{eq:c2.24}
\dfrac{\norm{X_1-x_{R_{m_2}}}}{\norm{X_1}}>\varepsilon,
\end{eqnarray}

is defined $ X_ {2} = x_{R_ {m_2}} $. Without loss of generality, it may be assumed that the second condition is fulfilled, then is taken $X_3= X_{R_ {m_3}} $ and are checked  the conditions \eqref{eq:c2.23} and \eqref{eq:c2.24} with the values $ X_1 $ and $ X_2 $. The  process described before is repeated for all values $ X_k=x_{R_{m_k}} $, with $ k \geq 4 $, and that generates a sequence $ \set{X_m}_ {m = 1} ^ t $, with $ t \leq r $, such that \label{pg:c2.03}

\begin{eqnarray*}
\dfrac{\norm{X_i-X_j}}{\norm{X_i}}>\varepsilon , & \forall i \neq j.
\end{eqnarray*}

\end{enumerate}

By following the steps described before to implement the fractional Newton method, a subset of  the solution set  of zeros, both real and complex, may be obtained from the function  $ f $. We will proceed to give an example where is found the solution set of zeros of a function  $f\in \nset{P}_n(\nset{R})$.

\begin{example}

Let the function:

\begin{eqnarray}\label{eq:0001}
\footnotesize
\begin{array}{rl}
f(x)=&  - 57.62x^{16} - 56.69x^{15}- 37.39x^{14} - 19.91x^{13} + 35.83x^{12}- 72.47x^{11} + 44.41x^{10} + 43.53x^{9}  \\
& + 59.93x^{8} - 42.9x^{7}   - 54.24x^{6} + 72.12x^{5}- 22.92x^{4}+ 56.39x^{3} + 15.8x^{2} +60.05x + 55.31,
\end{array}
\end{eqnarray}

then the following values are chosen to use the iteration function given by \eqref{eq:c2.20}

\begin{eqnarray*}
\footnotesize
\begin{array}{ccccc}
TOL=e-4,& L_{IT}=40, & \delta=0.5, &x_0=0.74, &  M=e+17,
\end{array}
\end{eqnarray*}

and using the fractional derivative given by \eqref{eq:c1.23}, we obtain the results of the Table \ref{tab:01}

\begin{table*}[!ht]
\centering
$
\footnotesize
\begin{array}{c|ccccc}
\toprule
&\alpha_m& {}^m\xi & \norm{{}^m \xi - {}^{m-1}\xi}_2 &\norm{f\left({}^m\xi \right)}_2& R_m \\ \midrule
1	&	-1.01346	&	 -1.3699527 	&	1.64700e-5	&	7.02720e-5	&	2	\\
2	&	 -0.80436	&	 -1.00133957 	&	9.82400e-5	&	4.36020e-5	&	2	\\
3	&	 -0.50138	&	 -0.62435277 	&	9.62700e-5	&	2.31843e-6	&	2	\\
4	&	  0.87611	&	  0.58999224 - 0.86699687i	&	3.32866e-7	&	6.48587e-6	&	11	\\
5	&	  0.87634	&	  0.36452488 - 0.83287821i	&	3.36341e-6	&	2.93179e-6	&	11	\\
6	&	  0.87658	&	 -0.28661369 - 0.80840642i	&	2.65228e-6	&	1.06485e-6	&	10	\\
7	&	  0.8943 	&	  0.88121183 + 0.4269622i 	&	1.94165e-7	&	6.46531e-6	&	14	\\
8	&	  0.89561	&	  0.88121183 - 0.4269622i 	&	2.87924e-7	&	6.46531e-6	&	11	\\
9	&	  0.95944	&	 -0.35983764 + 1.18135267i	&	2.82843e-8	&	2.53547e-5	&	24	\\
10	&	  1.05937	&	  1.03423976 	&	1.80000e-7	&	1.38685e-5	&	4	\\
11	&	  1.17776	&	 -0.70050491 - 0.78577099i	&	4.73814e-7	&	9.13799e-6	&	15	\\
12	&	  1.17796	&	 -0.35983764 - 1.18135267i	&	4.12311e-8	&	2.53547e-5	&	17	\\
13	&	  1.17863	&	 -0.70050491 + 0.78577099i	&	8.65332e-7	&	9.13799e-6	&	18	\\
14	&	  1.17916	&	  0.58999224 + 0.86699687i	&	7.05195e-7	&	6.48587e-6	&	12	\\
15	&	  1.17925	&	  0.36452488 + 0.83287821i	&	2.39437e-6	&	2.93179e-6	&	9	\\
16	&	  1.22278	&	 -0.28661369 + 0.80840642i	&	5.36985e-6	&	1.06485e-6	&	9	\\ \bottomrule
\end{array}
$
\caption{Results obtained using the iterative method \eqref{eq:c2.20}.}\label{tab:01}
\end{table*}

\end{example}

Although the fractional Newton method was originally defined for polynomials, the method can be extended to a broader class of functions, as shown in the following examples:

\begin{example}

Let the function:

\begin{eqnarray}\label{eq:c2.25}
\footnotesize
\begin{array}{c}
f(x)=  \sin(x)-\dfrac{3}{2x},
\end{array}
\end{eqnarray}

then the following values are chosen to use the iteration function given by \eqref{eq:c2.20}

\begin{eqnarray*}
\footnotesize
\begin{array}{ccccc}
TOL=e-4,& L_{IT}=40, & \delta=0.5, &x_0=0.26, &  M=e+6,
\end{array}
\end{eqnarray*}

and using the fractional derivative given by \eqref{eq:c1.23}, we obtain the results of Table \ref{tab:02}

\begin{table*}[!ht]
\centering
$
\footnotesize
\begin{array}{c|ccccc}
\toprule
&\alpha_m& {}^m\xi & \norm{{}^m \xi - {}^{m-1}\xi}_2 &\norm{f\left({}^m\xi \right)}_2& R_m \\ \midrule
1	&	-1.92915	&	1.50341195 	&	2.80000e-7	&	2.93485e-9	&	6	\\
2	&	 -0.07196	&	  -2.49727201 	&	9.99500e-5	&	6.53301e-9	&	  8	\\
3	&	 -0.03907	&	  -1.50341194 	&	6.29100e-5	&	4.37493e-9	&	  7	\\
4	&	  0.19786	&	 -18.92888307 	&	4.00000e-8	&	1.97203e-9	&	 20	\\
5	&	  0.20932	&	  -9.26211143 	&	9.60000e-7	&	4.77196e-9	&	 12	\\
6	&	  0.2097 	&	 -15.61173324 	&	5.49000e-6	&	2.05213e-9	&	 18	\\
7	&	  0.20986	&	  -12.6848988 	&	3.68000e-5	&	3.29282e-9	&	 15	\\
8	&	  0.21105	&	  -6.51548968 	&	9.67100e-5	&	2.05247e-9	&	 10	\\
9	&	  0.21383	&	 -21.92267274 	&	6.40000e-6	&	2.03986e-8	&	 24	\\
10	&	  1.19522	&	   6.51548968 	&	7.24900e-5	&	2.05247e-9	&	 13	\\
11	&	  1.19546	&	   9.26211143 	&	1.78200e-5	&	4.77196e-9	&	 14	\\
12	&	  1.19558	&	   12.6848988 	&	7.92100e-5	&	3.29282e-9	&	 14	\\
13	&	  1.19567	&	  15.61173324 	&	7.90000e-7	&	2.05213e-9	&	 12	\\
14	&	  1.1957 	&	  18.92888307 	&	1.00000e-8	&	1.97203e-9	&	 12	\\
15	&	  1.19572	&	  21.92267282 	&	1.46400e-5	&	5.91642e-8	&	 14	\\
16	&	  1.23944	&	     2.4972720 	&	6.30000e-7	&	9.43179e-10	&	 11		\\ \bottomrule
\end{array}
$
\caption{Results obtained using the iterative method \eqref{eq:c2.20}.}\label{tab:02}
\end{table*}

\end{example}

In the previous example, a subset of the solution set of zeros of the function \eqref{eq:c2.25} was obtained, because this function has an infinite amount of zeros. Using the fractional Newton method does not guarantee that all zeros of a function $ f $ can be found, leaving an initial condition $ x_0 $ fixed and varying the orders $ \alpha_m $ of the derivative. As in the classical Newton's method, finding most of the zeros of the function will depend on giving a proper initial condition $ x_0 $. If we want to find a larger subset of zeros of the function \eqref{eq:c2.25}, there are some strategies that are usually useful, for example:

\begin{enumerate}
\item To change the initial condition $ x_0 $.
\item To use a larger amount of $ \alpha_m $ values.
\item To increase the value of $ L_{IT} $.
\end{enumerate}

In general, the last strategy is usually the most useful, but this causes the  method \eqref{eq:c2.20} to become more costly, because it requires a longer runtime for all values $ \alpha_m $.

\begin{example}

Let the function:

\begin{eqnarray}
\footnotesize
\begin{array}{c}
f(x)= \left( x_1^2 + x_2^3 - 10, x_1^3 - x_2^2 - 1 \right)^T,
\end{array}
\end{eqnarray}

then the following values are chosen to use the iteration function given by \eqref{eq:c2.20}

\begin{eqnarray*}
\footnotesize
\begin{array}{ccccc}
TOL=e-4,& L_{IT}=40, & \delta=0.5, &x_0=(0.88,0.88)^T, &  M=e+6,
\end{array}
\end{eqnarray*}

and using the fractional derivative given by \eqref{eq:c1.23}, we obtain the results of Table \ref{tab:03} 

\begin{table*}[!ht]
\centering
$
\footnotesize
\begin{array}{c|cccccc}
&\alpha_m& {}^m\xi_1& {}^m\xi_2 &\norm{{}^m \xi - {}^{m-1} \xi }_2  &\norm{f\left({}^m\xi \right)}_2& R_m \\ \hline
1	&	-0.58345	&	0.22435853 + 1.69813926i	&	-1.13097646 + 2.05152306i	&	3.56931e-7	&	8.18915e-8	&	12	\\
2	&	 -0.50253	&	  0.22435853 - 1.69813926i	&	 -1.13097646 - 2.05152306i	&	1.56637e-6	&	8.18915e-8	&	 10	\\
3	&	  0.74229	&	 -1.42715874 + 0.56940338i	&	 -0.90233562 - 1.82561764i	&	1.13040e-6	&	7.01649e-8	&	 11	\\
4	&	  0.75149	&	  1.35750435 + 0.86070348i	&	  -1.1989996 - 1.71840823i	&	3.15278e-7	&	4.26428e-8	&	 12	\\
5	&	  0.76168	&	 -0.99362838 + 1.54146499i	&	   2.2675011 + 0.19910814i	&	8.27969e-5	&	1.05527e-7	&	 13	\\
6	&	  0.76213	&	 -0.99362838 - 1.54146499i	&	   2.2675011 - 0.19910815i	&	2.15870e-7	&	6.41725e-8	&	 15	\\
7	&	  0.77146	&	 -1.42715874 - 0.56940338i	&	 -0.90233562 + 1.82561764i	&	3.57132e-6	&	7.01649e-8	&	 15	\\
8	&	  0.78562	&	  1.35750435 - 0.86070348i	&	  -1.1989996 + 1.71840823i	&	3.16228e-8	&	4.26428e-8	&	 17	\\
9	&	  1.22739	&	  1.67784847 	&	  1.92962117 	&	9.99877e-5	&	2.71561e-8	&	  4	
\end{array}
$
\caption{Results obtained using the iterative method \eqref{eq:c2.20}.}\label{tab:03}
\end{table*}

\end{example}

\begin{example}

Let the function:

\begin{eqnarray}\label{eq:c2.29}
\footnotesize
\begin{array}{c}
f(x)= \left( x_1+x_2^2-37,x_1-x_2^2-5,x_1+x_2+x_3-3 \right) ^T,
\end{array}
\end{eqnarray}

then the following values are chosen to use the iteration function given by \eqref{eq:c2.20}

\begin{eqnarray*}
\footnotesize
\begin{array}{ccccc}
TOL=e-4,& L_{IT}=40, & \delta=0.5, &x_0=(4.35,4.35,4.35)^T, &  M=e+6,
\end{array}
\end{eqnarray*}

and using the fractional derivative given by \eqref{eq:c1.23}, we obtain the results of Table \ref{tab:05}

\begin{table*}[!ht]
\centering
$
\footnotesize
\begin{array}{c|ccccccc}
m&\alpha_m& {}^m\xi_1& {}^m\xi_2 & {}^m\xi_3 &\norm{{}^m \xi - {}^{m-1} \xi }_2 & \norm{f\left({}^m\xi \right)}_2& R_m \\ \hline
1	&	0.78928	&	-6.08553731 + 0.27357884i	&	0.04108101 + 3.32974848i	&	9.04445631 - 3.60332732i	&	6.42403e-5	&	3.67448e-8	&	14	\\
2	&	 0.79059	&	 -6.08553731 - 0.27357884i	&	  0.04108101 - 3.32974848i	&	  9.04445631 + 3.60332732i	&	1.05357e-7	&	3.67448e-8	&	 15	\\
3	&	 0.8166 	&	  6.17107462     	&	 -1.08216201     	&	 -2.08891261     	&	6.14760e-5	&	4.45820e-8	&	  9	\\
4	&	 0.83771	&	         6.0 	&	         1.0 	&	        -4.0 	&	3.38077e-6	&	0.0	&	  6	
\end{array}
$
\caption{Results obtained using the iterative method \eqref{eq:c2.20}.}\label{tab:05}
\end{table*}

\end{example}

\subsection{ Fractional Quasi-Newton Method}

Although the fractional Newton method is useful for finding multiple zeros of a  function $ f $, it has the disadvantage that in many cases calculating the fractional derivative of a function is not a simple task. To try to minimize this problem, we use that commonly for many definitions of the fractional derivative, the arbitrary order derivative of a constant is not always zero, that is,

\begin{eqnarray}\label{eq:c2.30}
\footnotesize
\begin{array}{cc}
\der{\partial}{x_j}{\alpha}c\neq 0, & c=constant.
\end{array}
\end{eqnarray}

Then, we may define the function

\begin{eqnarray}\label{eq:c2.31}
g(x)=f(x_0)+f^{(1)}(x_0)x,
\end{eqnarray}

it should be noted that the previous function is almost a linear approximation of the function $ f $ in the initial condition $ x_0 $. Then, for any fractional derivative that satisfies the condition \eqref{eq:c2.30}, and using \eqref{eq:c2.26} the \textbf {Fractional Quasi-Newton  Method} may be defined as

\begin{eqnarray}\label{eq:c2.32}
\begin{array}{cc}
x_{i+1}:=\Phi(\alpha,x_i)=x_i-\left(  Q_{g,\beta}(x_i)  \right)^{-1}f(x_i), & i=0,1,2,\cdots.,
\end{array}
\end{eqnarray}

with $ Q_ {g, \beta} (x_i) $ given by the following matrix 

\begin{eqnarray}
Q_{g,\beta}(x_i):=\left(
\partial_{j}^{\beta( \alpha,(x_i)_j )} g_k(x_i)\right),
\end{eqnarray}

where $(x_i)_j$ is the component $j$-th of the value $x_i$, and $ \beta(\alpha,(x_i)_j) $ is defined as follows

\begin{eqnarray}\label{eq:c2.34}
\normalsize
\beta(\alpha,( x_i)_j):=\left\{
\begin{array}{cc}
\alpha, & \mbox{if } \abs{(x_i)_j}\neq 0 ,\\
1, & \mbox{if } \abs{(x_i)_j}=0.
\end{array}\right.
\end{eqnarray}

Since the iteration function \eqref {eq:c2.32} does not satisfy the condition \eqref{eq:c2.17}, then any sequence $ \set{x_i} _ {i = 0} ^ \infty $ generated by this iteration function has at most one order of convergence (at least) linear. As a consequence, the speed of convergence is slower compared to what would be obtained when using \eqref{eq:c2.20}, and then it is necessary to use a larger value of $ L_ {IT} $. It should be mentioned that the value $ \alpha = 1 $ in \eqref{eq:c2.34}, is not taken to try to guarantee an order of convergence as in \eqref{eq:c2.21}, but to avoid the discontinuity that is generated when using the fractional derivative of constants in the value $ x = 0 $. An example is given using the fractional quasi-Newton  method, where is found a subset of the  solution set of zeros  of the function $ f $.

\begin{example}

Let the function:

\begin{eqnarray}\label{eq:c2.33}
\footnotesize
\begin{array}{c}
f(x)= \left(\dfrac{1}{2}\sin(x_1x_2)-\dfrac{x_2}{4\pi}-\dfrac{x_1}{2}, \left( 1-\dfrac{1}{4\pi}\right)\left(e^{2x_1}-e \right)+\dfrac{e}{\pi}x_2-2ex_1 \right)^T,
\end{array}
\end{eqnarray}

then the following values are chosen to use the iteration function given by \eqref{eq:c2.32}

\begin{eqnarray*}
\footnotesize
\begin{array}{cccc}
TOL=e-4,& L_{IT}=200, &x_0=(1.52,1.52)^T, &  M=e+6,
\end{array}
\end{eqnarray*}

and using the fractional derivative given by \eqref{eq:c1.23}, we obtain the results of Table \ref{tab:06}

\begin{table*}[!ht]
\centering
$
\footnotesize
\begin{array}{c|cccccc}
&\alpha_m& {}^m\xi_1& {}^m\xi_2 &\norm{{}^m \xi - {}^{m-1} \xi }_2  &\norm{f\left({}^m\xi \right)}_2& R_m \\ \hline
1	&	-0.28866	&	2.21216549 - 13.25899819i	&	0.41342314 + 3.94559327i	&	1.18743e-7	&	9.66251e-5	&	163	\\
2	&	  1.08888	&	  1.29436489      	&	  -3.13720898     	&	1.89011e-6	&	9.38884e-5	&	  51	\\
3	&	  1.14618	&	  1.43395246      	&	  -6.82075021     	&	2.24758e-6	&	9.74642e-5	&	  94	\\
4	&	  1.33394	&	  0.50000669      	&	   3.14148062     	&	9.74727e-6	&	9.99871e-5	&	 133	\\
5	&	  1.35597	&	  0.29944016      	&	   2.83696105     	&	8.55893e-5	&	4.66965e-5	&	   8	\\
6	&	  1.3621 	&	   1.5305078      	&	 -10.20223066     	&	2.38437e-6	&	9.88681e-5	&	 120	\\
7	&	  1.37936	&	  1.60457254      	&	 -13.36288413     	&	2.32459e-6	&	9.52348e-5	&	  93	\\
8	&	  1.88748	&	 -0.26061324 	&	   0.62257513 	&	2.69146e-5	&	9.90792e-5	&	  21	
\end{array}
$
\caption{Results obtained using the iterative method \eqref{eq:c2.32}.}\label{tab:06}
\end{table*}

\end{example}

It is worth mentioning that the function \eqref{eq:c2.33} has infinite zeros, so it will be used with the following method.

\subsection{ Fractional Pseudo-Newton Method}

A conventional way in which Newton's method is usually deduced is based on a linear approximation of a  function $ f $ in a value $ x_i $. To perform the aforementioned, we may consider the Taylor polynomial of the function  $ f $ around a value $ x_i $, and disregarding higher order terms we have that

\begin{eqnarray}\label{eq:c2.36}
f(x)\approx f(x_i)+f^{(1)}(x_i)(x-x_i),
\end{eqnarray}

then, assuming that $ \xi $ is a zero of $ f $, from the previous expression we have that

\begin{eqnarray*}
\xi \approx x_i- \dfrac{f(x_i)}{f^{(1)}(x_i)},
\end{eqnarray*}

as consequence, a sequence $ \set{x_i}_{i = 0} ^ \infty $ that approximates the value $ \xi $ may be generated using the iteration function

\begin{eqnarray*}
x_{i+1}:=\Phi(x_i)=x_i-\dfrac{f(x_i)}{f^{(1)}(x_i)}, & i=0,1,2,\cdots.
\end{eqnarray*}

However, the equation \eqref{eq:c2.36} is not the only way to generate a linear approximation to the  function $ f $ in the point $ x_i $, in general it may be taken as

\begin{eqnarray}\label{eq:c2.37}
f(x)\approx f(x_i)+m(x-x_i),
\end{eqnarray}

where $ m $ is any constant value of a slope, that allows the approximation  \eqref{eq:c2.37} to the  function $ f $ to be valid. The previous equation allows to obtain the following iteration function

\begin{eqnarray}\label{eq:c2.38}
x_{i+1}:=\Phi(x_i)=x_i- \dfrac{f(x_i)}{m},&  i=0,1,2,\cdots,
\end{eqnarray}

which originates the \textbf{Parallel Chord Method} \cite{ortega1970iterative}. The iteration function \eqref{eq:c2.38} can be generalized to larger dimensions as follows

\begin{eqnarray}\label{eq:c2.39}
x_{i+1}=\Phi(x_i)=x_i-\left(\dfrac{1}{m}I_n\right) f(x_i),&  i=0,1,2,\cdots,
\end{eqnarray}

where $ I_n $, is the identity matrix of $ n \times n $. The previous equation implies that it is enough to apply the method \eqref{eq:c2.38} component to component for the case in several variables. 

Using as a basis the idea of \eqref{eq:c2.39}, and considering any fractional derivative that satisfies the condition \eqref{eq:c2.30}, we can define the \textbf{Fractional Pseudo-Newton Method} as follows

\begin{eqnarray}\label{eq:c2.40}
\begin{array}{cc}
x_{i+1}:=\Phi\left(\alpha,x_i \right)=x_i- P_{\epsilon,\beta}(x_i)f(x_i), & i=0,1,2\cdots,
\end{array} 
\end{eqnarray}

with $ P_ {\epsilon, \beta} (x_i) $ given by the following matrix evaluated at the value $ x_i $

\begin{eqnarray}
P_{\epsilon,\beta}(x_i):=\left(
\partial_{j}^{\beta( \alpha,(x_i)_j )} \delta_{kj} + \epsilon \delta_{kj}\right)_{x_i},
\end{eqnarray}

where $(x_i)_j$ is the component $j$-th of the value $x_i$, with $\delta_{kj}$ the Kronecker delta, $\epsilon$ a positive constant $\ll 1$, and  $ \beta (\alpha,(x_i)_j) $ is given by \eqref{eq:c2.34}.

Two examples are given using the fractional pseudo-Newton method, where is found a subset of the solution set of zeros of some functions. 
 
\newpage

\begin{example}

Let the function:

\begin{eqnarray*}
\footnotesize
\begin{array}{c}
f(x)= \left(\dfrac{1}{2}\sin(x_1x_2)-\dfrac{x_2}{4\pi}-\dfrac{x_1}{2}, \left( 1-\dfrac{1}{4\pi}\right)\left(e^{2x_1}-e \right)+\dfrac{e}{\pi}x_2-2ex_1 \right)^T,
\end{array}
\end{eqnarray*}

then the following values are chosen to use the iteration function given by \eqref{eq:c2.40}

\begin{eqnarray*}
\footnotesize
\begin{array}{ccccc}
TOL=e-4,& L_{IT}=200,  &x_0=(1.03,1.03)^T, &  M=+6,&\epsilon=e-3,
\end{array}
\end{eqnarray*}

and using the fractional derivative given by \eqref{eq:c1.23}, we obtain the results of Table \ref{tab:07} 

\begin{table*}[!ht]
\centering
$
\footnotesize
\begin{array}{c|cccccc}
\toprule
&\alpha_m& {}^m\xi_1& {}^m\xi_2 &\norm{{}^m \xi - {}^{m-1} \xi }_2  &\norm{f\left({}^m\xi \right)}_2& R_m \\ 
\midrule
1	&	0.78562	&	 1.03499277 - 0.53982128i	&	5.41860852 + 4.04164098i	&	5.62354e-6	&	8.38442e-5	&	66	\\
2	&	 0.78987	&	  0.29945564        	&	  2.83683317        	&	1.09600e-5	&	9.63537e-5	&	  88	\\
3	&	 0.82596	&	 -0.26054499        	&	  0.62286899        	&	5.66073e-5	&	9.87374e-5	&	 140	\\
4	&	 0.82671	&	  -0.1561964 - 1.02056003i	&	  2.26280132 - 5.71855964i	&	4.32875e-6	&	9.51178e-5	&	 194	\\
5	&	 0.83158	&	  1.03499697 + 0.53981525i	&	  5.41862187 - 4.04161017i	&	3.94775e-6	&	8.80344e-5	&	  84	\\
6	&	 0.85861	&	  1.16151359 - 0.69659512i	&	  8.27130854 + 6.3096935i 	&	2.14707e-6	&	9.38721e-5	&	 164	\\
7	&	 1.15911	&	  1.48131686        	&	 -8.38362876        	&	1.20669e-6	&	9.56674e-5	&	 191	\\
8	&	 1.24977	&	 -1.10844524 + 0.10906317i	&	 -4.18608959 + 0.66029327i	&	3.71508e-6	&	9.81146e-5	&	 164	\\
9	&	 1.25662	&	 -1.10844605 - 0.10906368i	&	 -4.18608629 - 0.66029181i	&	3.69483e-6	&	9.66271e-5	&	 170	\\
10	&	 1.26128	&	  1.33741853  	&	 -4.14026671  	&	1.89913e-5	&	8.51053e-5	&	  67	\\
\bottomrule		
\end{array}
$
\caption{Results obtained using the iterative method \eqref{eq:c2.40}}\label{tab:07}
\end{table*}

\end{example}

\begin{example}

Let the function:

\begin{eqnarray}\label{eq:c2.29}
\footnotesize
\begin{array}{c}
f(x)= \left(  -3.6 x_3\left( x_1^3x_2 +1 \right) - 3.6\cos\left( x_2^2 \right) + 10.8, 
-1.6x_1\left( x_1 +x_2^3x_3\right) - 1.6\sinh\left( x_3\right) + 6.4, 
-4x_2\left( x_1x_3^3 +1 \right) - 4\cosh\left( x_1 \right) + 24 \right)^T,
\end{array}
\end{eqnarray}

then the following values are chosen to use the iteration function given by \eqref{eq:c2.40}

\begin{eqnarray*}
\footnotesize
\begin{array}{ccccc}
TOL=e-4, & L_{IT}=200, &x_0=(1.12,1.12,1.12)^T, &  M=e+6 ,& \epsilon=e-3,
\end{array}
\end{eqnarray*}

and using the fractional derivative given by \eqref{eq:c1.23}, we obtain the results of Table \ref{tab:08}

\begin{table*}[!ht]
\centering
$
\footnotesize
\begin{array}{c|ccccccc}
\toprule
&\alpha_m& {}^m\xi_1& {}^m\xi_2 &{}^m\xi_3&\norm{{}^m \xi - {}^{m-1} \xi }_2  &\norm{f\left({}^m\xi \right)}_2& R_m \\ 
\midrule
1	&	0.96743	&	0.38147704 + 1.10471108i	&	-0.43686196 - 1.3473184i 	&	-0.38512615 - 1.4903386i 	&	2.41936e-6	&	7.07364e-5	&	57	\\
2	&	 0.96745	&	 -0.78311553 + 0.96791081i	&	 -0.58263802 + 1.2592471i 	&	  0.18175185 - 1.49135484i	&	3.85644e-6	&	8.58385e-5	&	 37	\\
3	&	 0.96766	&	  0.71500126 - 1.02632085i	&	  0.53575431 - 1.314774i  	&	  0.45273307 - 1.35710557i	&	3.01511e-6	&	8.34643e-5	&	 41	\\
4	&	 0.9677 	&	 -0.34118928 + 1.19432023i	&	  0.37199268 - 1.40125985i	&	 -0.63215137 + 1.3074313i 	&	2.72698e-6	&	8.69377e-5	&	 49	\\
5	&	 0.96796	&	  0.71500155 + 1.0263218i 	&	  0.53575489 + 1.31477495i	&	  0.45273303 + 1.35710453i	&	2.52069e-6	&	7.11216e-5	&	 34	\\
6	&	 0.97142	&	 -0.34118945 - 1.19432007i	&	  0.37199303 + 1.40125973i	&	 -0.63215109 - 1.3074314i 	&	2.32465e-6	&	8.66652e-5	&	 61	\\
7	&	 0.9718 	&	  0.38147878 - 1.10471296i	&	 -0.43686073 + 1.34732029i	&	  -0.3851262 + 1.49033466i	&	2.38466e-6	&	8.53472e-5	&	 50	\\
8	&	 0.97365	&	 -0.78311508 - 0.96791138i	&	 -0.58263753 - 1.2592475i 	&	  0.18175161 + 1.49135503i	&	1.99078e-6	&	7.57517e-5	&	 51	\\
9	&	 1.03148	&	  1.34508926      	&	 -1.29220278      	&	 -1.44485467      	&	3.68616e-6	&	9.58451e-5	&	 59	\\
10	&	 1.04155	&	 -1.43241693 	&	  1.27535274  	&	 -1.11183615 	&	4.06891e-6	&	8.95830e-5	&	 48	\\
\bottomrule
\end{array}
$
\caption{Results obtained using the iterative method \eqref{eq:c2.40}}\label{tab:08}
\end{table*}

\end{example}

\section{Conclusions}

The fractional Newton-Raphson  method and its variants are useful for finding multiple solutions of nonlinear systems, in the complex space using real initial conditions. However, it should be clarified that they present some advantages and disadvantages between each of them, for example, although the fractional Newton method generally has an order of convergence (at least) quadratic, this method has the disadvantage that it is not an easy task to find the fractional Jacobian matrix for many functions, and also the need to reverse this matrix must be added for each new iteration. But it has an advantage over the other methods, and this is because it can be used with few iterations, which allows occupying a greater number of $ \alpha_m $ values belonging to the partition of the interval $ (- 2,2) $. 

The quasi-Newton fractional method has the advantage that the fractional Jacobian matrix with which it works, compared to the fractional Newton method, is easy to obtain. But a disadvantage is that the method may have at most an order of convergence (at least) linear, so the speed of convergence is lower and it is necessary to use a greater number of iterations to ensure success in the search for solutions. As a consequence, the method is more costly because it requires a longer runtime to use all values $ \alpha_m $. An additional advantage of the method is that if the initial condition is close enough to a solution, its behavior is very similar to fractional Newton method and may converge with a relatively small number of iterations, but it still has the disadvantage that we need to reverse a matrix in each iteration.

Finally, the fractional pseudo-Newton method seems to solve the problem of the need to invert a matrix in each iteration that are present in the previous methods. However, this method, as the fractional quasi-Newton method, may  have at most an order of convergence (at least) linear, and it has a lower speed of convergence  compared to the latter. As a consequence, it is necessary to use a larger number of values $ \alpha_m $ and a greater number of iterations, then we require a longer runtime to use all values to find solutions, so it can be considered as a slow and costly method. But its advantage over the other methods is that it is not necessary to calculate and invert any fractional Jacobian matrix. 

The methods may solve some nonlinear systems and are really efficient to find multiple solutions, both real and complex, using real initial conditions. It should be mentioned that these methods are extremely recommended in systems that have infinite solutions or a large number of them.

\bibliography{Biblio}
\bibliographystyle{unsrt}
\nocite{brambila2019proposal}
\nocite{martinez2017applications1}
\nocite{martinez2017application2}

\end{document}